\documentclass[aop,preprint]{imsart}

\usepackage{enumerate}
\usepackage{amssymb}
\usepackage{amsthm,amsmath}
\usepackage{amsfonts}
\usepackage{mathrsfs}
\usepackage{graphicx}
\usepackage{mathrsfs}
\usepackage{stmaryrd}
\usepackage{bbm}
\usepackage{cases}
\usepackage{amsfonts}
\usepackage{txfonts}

\arxiv{arXiv:1407.5834}

\startlocaldefs

\newcommand{\be}{\begin{eqnarray}}
\newcommand{\ee}{\end{eqnarray}}
\newcommand{\ce}{\begin{eqnarray*}}
\newcommand{\de}{\end{eqnarray*}}
\newtheorem{theorem}{Theorem}[section]
\newtheorem{lemma}[theorem]{Lemma}
\newtheorem{remark}[theorem]{Remark}
\newtheorem{definition}[theorem]{Definition}
\newtheorem{proposition}[theorem]{Proposition}
\newtheorem{Examples}[theorem]{Example}
\newtheorem{corollary}[theorem]{Corollary}

\def\e{{\mathrm{e}}}
\def\eps{\varepsilon}
\def\t{\tau}

\def\p{\partial}

\def\[{{\Big[}}
\def\]{{\Big]}}
\def\<{{\langle}}
\def\>{{\rangle}}
\def\({{\Big(}}
\def\){{\Big)}}

\def\bx{{\mathbf{x}}}

\def\dif{{\mathord{{\rm d}}}}
\def\dis{{\mathord{{\rm \bf d}}}}

\def\no{\nonumber}
\def\={&\!\!=\!\!&}
\def\bt{\begin{theorem}}
\def\et{\end{theorem}}
\def\bl{\begin{lemma}}
\def\el{\end{lemma}}
\def\br{\begin{remark}}
\def\er{\end{remark}}

\def\bd{\begin{definition}}
\def\ed{\end{definition}}
\def\bp{\begin{proposition}}
\def\ep{\end{proposition}}
\def\bc{\begin{corollary}}
\def\ec{\end{corollary}}
\def\bx{\begin{Examples}}
\def\ex{\end{Examples}}

\def\cM{{\mathcal M}}

\def\mA{{\mathbb A}}

\def\mE{{\mathbb E}}

\def\mL{{\mathbb L}}

\def\mN{{\mathbb N}}

\def\mP{{\mathbb P}}
\def\mQ{{\mathbb Q}}
\def\mR{{\mathbb R}}

\def\mT{{\mathbb T}}

\def\mV{{\mathbb V}}
\def\mW{{\mathbb W}}

\def\sE{{\mathscr E}}
\def\sF{{\mathscr F}}

\def\sL{{\mathscr L}}

\def\geq{\geqslant}
\def\leq{\leqslant}

\def\t{\mathord{{\rm tr}}}

\endlocaldefs

\begin{document}

\begin{frontmatter}

\title{Sobolev differentiable flows of SDEs with local Sobolev and super-linear growth coefficients}
\runtitle{Sobolev differentiable flows of SDEs}

\author{\fnms{Longjie} \snm{Xie}\corref{}\ead[label=e1]{xlj.98@whu.edu.cn}}
\address{School of Mathematics and Statistics,\\ Wuhan University,\\ Wuhan, Hubei 430072, P.R.China\\ \printead{e1}}
\and
\author{\fnms{Xicheng} \snm{Zhang}\ead[label=e2]{XichengZhang@gmail.com}}
\address{School of Mathematics and Statistics,\\ Wuhan University,\\ Wuhan, Hubei 430072, P.R.China\\ \printead{e2}}
\affiliation{Wuhan University}

\runauthor{}

\begin{abstract}
By establishing a characterization for Sobolev differentiability of random fields,
we prove the weak differentiability of solutions to stochastic differential equations with local Sobolev and super-linear growth coefficients
with respect to the starting point. Moreover, we also study the strong Feller property and the irreducibility to the associated diffusion semigroup.
\end{abstract}



\end{frontmatter}

\section{Introduction and Main Results}

Consider the following stochastic differential equation (SDE) in $\mR^d$:
\begin{align}
\dif X_t(x) =b(t,X_t(x))\dif t+\sigma(t,X_t(x))\dif W_t, \quad X_0(x)=x,\label{1}
\end{align}
where $\sigma:\mR_+\times\mR^d\rightarrow\mR^d\otimes\mR^m$ and $b:\mR_+\times\mR^d\rightarrow\mR^d$ are two measurable functions, $(W_t)_{t\geq 0}$ is
an $m$-dimensional standard Brownian motion defined on some probability space $(\Omega,\sF,\mP)$.
It is a classical result that if the coefficients are global Lipschitz continuous and have linear growth in $x$ uniformly with respect to $t$,
then SDE (\ref{1}) admits a unique global strong solution which forms a stochastic flow of homeomorphisms on $\mR^d$ (cf. \cite{Ku}).
However, in many applications, the Lipschitz continuity and linear growth condition imposed on the coefficients are broken
(see \cite{Di-Li}, \cite{F-G-P}, \cite{Kr-Ro}, \cite{Pr-Ro} and references therein).
Notice that in the deterministic case (i.e., $\sigma\equiv0$), SDE (\ref{1}) becomes an ordinary differential equation (ODE):
\begin{align}
x'(t)=b\big(t,x(t)\big),\quad x(0)=x_0.    \label{ode}
\end{align}
A unique regular Lagrangian flow was constructed in \cite{Di-Li} by DiPerna and Lions for Sobolev vector fields with bounded divergence (see also \cite{Cr-De}
for a direct argument).
This result was later extended by Ambrosio in \cite{Am} to BV vector fields with bounded divergence (see also \cite{Fa-Lu-Th}, \cite{Zh0} and \cite{Ch-Ja}
for stochastic extensions). It is emphasized that the solvability of \eqref{ode} in the DiPerna-Lions theory is only for Lebesgue almost all starting point $x_0$.
An interesting phenomena is that when $\sigma\neq0$ is nondegenerate, the noise term will play some regularization effect and SDE (\ref{1}) can be well-posed
for quite singular drift $b$ and for every starting point $x$.
In the past decades, there is increasing interest in the study of the strong solutions and their properties to SDEs (\ref{1}) with irregular coefficients.
Let us briefly recall some well-known results in this direction.

\vspace{2mm}

In the additive noise case (i.e., $\sigma_t(x)=\sigma_t$ is non-degenerate), when $b$ is bounded and measurable, Veretennikov \cite{Ve}
first proved that SDE (\ref{1}) has a unique global strong solution $X_t(x)$. Recently,
it was shown in \cite{Mo-Ni-Pr} that $X_t(\cdot)$ lies in the space
$\cap_{p\geq 1}L^2\big(\Omega; W^{1,p}_\rho(\mR^d)\big)$, where $W^{1,p}_\rho(\mR^d)$
denotes weighted Sobolev space with weight $\rho$ possessing finite $p$-th moment with respect to the Lebesgue measure in $\mR^d$.
When $b\in L^q_{loc}(\mR_+;L^p(\mR^d))$ for some $p,q\in(1,\infty)$ with
$\frac{d}{p}+\frac{2}{q}<1,$ using estimates of solutions to the associated PDE,
the existence and uniqueness of a global strong solution $X_t(x)$ for SDE (\ref{1}) were obtained by Krylov and R\"ockner in \cite{Kr-Ro}.
Under the same condition, Fedrizzi and Flandoli \cite{Fe-Fl-1,Fe-Fl-2, Fe-Fl-3}
proved that the map $x\rightarrow X_t(x)$ is $\alpha$-H\"older continuous for any $\alpha\in(0,1)$, and is also Sobolev differentiable.
We also mention that the Sobolev regularity of the strong solution in spacial variable enables us to study the associated stochastic transport equation
since it is closely related to SDE (\ref{1}) by the inverse flow of the strong solution, see \cite{Fe-Fl-3,F-G-P,Mo-Ni-Pr, Re}
and references therein. Moreover, Bismut-Elworthy-Li's formula was also established in \cite{M-M-N-P-Z}
by using the Sobolev and Malliavin differentiabilities of strong solutions with respect to the initial values and sample paths respectively.

\vspace{2mm}

In the multiplicative noise case: if the SDE is time homogeneous, supposing that $\sigma(x)$, $b(x)$ are in $C^2(\mR^d)$, and $\nabla\sigma$
and $\nabla b$ have some mild growth at infinity, by investigating the corresponding derivative flow equation, Li  \cite{Li} studied the strong completeness for SDE (\ref{1}),
i.e., $(t,x)\mapsto X_t(x)$ admits a bicontinuous version. More recently, this result was extended to the case of Sobolev coefficients in \cite{Ch-Li}
and the Sobolev regularity of solutions with respect to the initial value was also studied. The main argument in \cite{Ch-Li} is the mollifying approximation to SDE (\ref{1})
and the key point is to prove some uniform estimates of the solution to the corresponding derivative flow equation.
It is emphasized that in  \cite{Li, Ch-Li}, $\nabla\sigma$ and $\nabla b$ are not necessarily bounded.
Very recently, Zhang \cite{Zh3, Zh1, Zh4} proved under the assumptions that $\sigma$ is bounded, uniformly elliptic and uniformly continuous in $x$
locally uniformly with respect to $t$, and $|b|, |\nabla\sigma|\in L^q_{loc}(\mR_+\times\mR^d)$ for $q>d+2$,
there exists a unique strong solution $X_t(x)$ to (\ref{1}) up to the explosion time $\zeta(x)$ for every $x\in\mR^d$. Meanwhile,
under the {\it global} assumptions that $|b|, |\nabla\sigma|\in L^q_{loc}(\mR_+; L^q(\mR^d))$ for $q>d+2$,
the solution $\{X_t(x)\}$ forms a stochastic flow of homeomorphisms on $\mR^d$, and $x\mapsto X_t(x)$ is Sobolev differentiable.

\vspace{2mm}

In this paper, we will establish the Sobolev regularity of strong solutions with respect to the initial value, as well as the strong Feller property and irreducibility,
to SDE (\ref{1}) with some {\it local} Sobolev and super-linear growth coefficients. For this purpose,
we first establish a useful characterization for Sobolev differentiability of random fields in terms of their moment estimates, which has independent interest.

\bt\label{th1}
Let $U\subset\mR^d$ be a bounded $C^1$-domain and $f\in L^{q}\big(U;L^p(\Omega; L^r(T))\big)$ for some $p\in(1,\infty)$ and $q,r\in(1,\infty]$.
Then $f\in\mW^{1,q}\big(U;L^p(\Omega; L^r(T))\big)$ (see (\ref{w}) below for a definition)
if and only if there exists a nonnegative measurable function $g\in L^q(U)$ such that for Lebesgue-almost all $x,y\in U$,
\begin{align}
\|f(x,\cdot)-f(y,\cdot)\|_{L^p(\Omega; L^r(T))}\leq |x-y|\big(g(x)+g(y)\big).\label{Mo23}
\end{align}
Moreover, if (\ref{Mo23}) holds, then for Lebesgue-almost all $x\in U$,
\begin{align*}
\|\p_if(x,\cdot)\|_{L^p(\Omega; L^r(T))}\leq 2g(x),\ \ i=1,\cdots,d,
\end{align*}
where $\p_i f$ is the weak partial derivative of $f$ with respect to the  $i$-th spacial variable.
\et

The advantage of this characterization lies in that, when we want to show the Sobolev regularity of
the strong solution $X_t(x)$ to SDE (\ref{1}) with respect to $x$,
we just need to estimate the $p$-th moment of $X_t(x)-X_t(y)$, which is much easier to be handled for SDEs
(see a recent work \cite{Wa-Zh} for an application of the above characterization).
It should be noticed that in all previous works (see \cite{Ch-Li, Fe-Fl-3, Mo-Ni-Pr}), the argument of mollifying coefficients is used to obtain the Sobolev differentiability of
strong solutions. This usually leads to some complicated limiting procedures.
Here an interesting open question is that whether we can extend the above characterization to the infinite dimensional case in somehow so that
it can be used to the SDE in Hilbert spaces as studied in \cite{AB1, AB2, AB3, AB4}.

\vspace{2mm}

Now, we turn to the study of SDE \eqref{1} and make the following assumptions on $\sigma$ and $b$.
\vspace{2mm}
\begin{enumerate}[{\bf (H1)}]
\item {\bf (Local Sobolev integrability)} $\sigma$ is locally uniformly continuous in $x$ and locally uniformly with respect to $t\in\mR_+$,
and for some $q>d+2$,
$$
b\in L^q_{loc}(\mR_+\times\mR^d),\,\,\nabla\sigma\in L^q_{loc}(\mR_+\times\mR^d),
$$
and for some $C_1,\gamma_1>0$, $\alpha'\in(0,\alpha)$ and for all $t\geq 0$, $x\in\mR^d$, $\xi\in\mR^m$,
\begin{align}
|\sigma(t,x)\xi|\geq |\xi|\Big(1_{\alpha>0}\exp\big\{-C_1(1+|x|^2)^{\alpha'}\big\}+1_{\alpha=0}C_1(1+|x|^2)^{-\gamma_1}\Big),\label{c2}
\end{align}
where $\alpha$ is the same as in (\ref{c1}) below.
\item {\bf (Super-linear growth)} For some $\alpha\in[0,1]$ and for all $\kappa>0$,  there exist a constant $C_\kappa\in\mR$
and a nonnegative function $F_{\kappa}(t,x)\in L^{q'}_{loc}(\mR_+\times\mR^d)$ with some $q'>d+1$ such that for all $t\geq 0$ and $x,y\in\mR^d$,
\begin{align}
\<x,b(t,x)\>+\kappa\big(1+|x|^{2}\big)^{\alpha}\|\sigma(t,x)\|^2\leq C_\kappa\big(1+|x|^{2}\big), \label{c1}
\end{align}
\begin{align} \label{c4}
\begin{split}
&\<x-y,b(t,x)-b(t,y)\>+\kappa\|\sigma(t,x)-\sigma(t,y)\|^2\\
&\qquad\quad\leq |x-y|^2\Big(F_{\kappa}(t,x)+F_{\kappa}(t,y)\Big), 
\end{split}
\end{align}
and there exist $\alpha'\in[0,\alpha)$, $R_0>0$ and $C_2, \gamma_2, C_3>0$ such that for all $t\geq 0$ and $x\in\mR^d$,
\begin{align}\label{c3} 
|b(t,x)|+\|\sigma(t,x)\|\leq 1_{\alpha>0}\exp\big\{C_2(1+|x|^2)^{\alpha'}\big\}+1_{\alpha=0}C_2(1+|x|^2)^{\gamma_2},
\end{align}
and for all $t\geq 0$ and $|x|\geq R_0$,
\begin{align}
F_{\kappa}(t,x)\leq C_3\Big(1_{\alpha>0}(1+|x|^{2})^{\alpha'}+1_{\alpha=0}\log(1+|x|^2)\Big).    \label{cf}
\end{align}
\end{enumerate}

Here neither uniformly elliptic nor the global $L^q$-integrability conditions are assumed on $b$ and $\sigma$.
Notice that  $q>d+2$ in {\bf (H1)} is almost optimal due to Krylov and R\"ockner's sharp condition $\frac{d}{p}+\frac{2}{q}<1$.
Our first main result of this paper is

\bt\label{main1}
Under ${\bf{(H1)}}$ and ${\bf{(H2)}}$, there exists a unique global strong solution $X_t(x)$ to SDE (\ref{1}) so that
$(t,x)\mapsto X_t(x)$ is continuous. Moreover, we have the following conclusions:
\begin{enumerate}[\bf (A)]
\item For each $t>0$ and almost all $\omega$, the mapping $x\mapsto X_t(x,\omega)$ is Sobolev differentiable,
and for any $T>0$ and $p\geq 1$, there are constants $C,\gamma>0$ such that for Lebesgue-almost all $x\in\mR^d$,
\begin{align}
\sup_{t\in[0,T]}\mE|\nabla X_t(x)|^p\leq C\Big(1_{\alpha>0}\e^{(1+|x|^2)^\alpha}+1_{\alpha=0}(1+|x|^2)^\gamma\Big),   \label{so11}
\end{align}
where $\nabla$ denotes the gradient in the distributional sense, and $\alpha$ is the same as in (\ref{c1}).
\item If in addition, we assume that for some $F_0(t,x)\in L^{q'}_{loc}(\mR_+\times\mR^d)$ with $q'>d+1$,
\begin{align}\label{Sigma}
\|\sigma(t,x)-\sigma(t,y)\|^2\leq |x-y|^2\Big(F_{0}(t,x)+F_{0}(t,y)\Big),
\end{align}
where $F_0$ also satisfies \eqref{cf}, then \eqref{so11} can be strengthened as
\begin{align}
\mE\left(\mathrm{ess.}\sup_{t\in[0,T]}|\nabla X_t(x)|^p\right)\leq C\Big(1_{\alpha>0}\e^{(1+|x|^2)^\alpha}+1_{\alpha=0}(1+|x|^2)^\gamma\Big). \label{so1}
\end{align}
\item For each $t>0$ and any bounded measurable function $f$ on $\mR^d$,
$$
x\mapsto \mE f(X_t(x))\mbox{ is continuous}.
$$
\item For each open set $A\subset\mR^d$ and $t>0$, $x\in\mR^d$,
$$
\mP\{\omega: X_t(x,\omega)\in A\}>0.
$$
\end{enumerate}
\et
\br
If $b$ and $\sigma$ are time independent, then the above {\bf (C)} means that the semigroup defined by $P_tf(x):=\mE f(X_t(x))$ is strong Feller, and
the above {\bf (D)} means that $P_t$ is irreducible. In particular, {\bf (C)} and {\bf (D)}
imply the uniqueness of the invariant measures associated to $(P_t)_{t\geq 0}$ (if it exists). See \cite{Ce, Pr-Ro} for applications.
\er
\br
Assumptions (\ref{c1}) and (\ref{c4}) are classical coercivity and monotonicity conditions when $\kappa=\frac{1}{2}$, $\alpha=0$
and $F_\kappa(t,x)=constant$ in (\ref{c1}) and (\ref{c4}).
In this case, if, in addition that $b$ and $\sigma$ are continuous in $x$, then the existence and uniqueness of strong solutions to SDE (\ref{1}) are classical (cf. \cite{Pr-Ro}).
However, under non degenerated assumption (\ref{c2}), we can drop the continuity assumption on drift $b$. Moreover, our estimate \eqref{so1}
is stronger than the well-known results (cf. \cite{Ce, Fe-Fl-3, Mo-Ni-Pr, Ch-Li}) since the essential supremum norm with respect to the time variable is taken in the expectation.
\er

It is well-known that under ${\bf{(H1)}}$, SDE (\ref{1}) admits a unique local strong solution.
We will show in Lemma \ref{Le31} below that SDE (\ref{1}) with coefficients satisfying (\ref{c1}) does not explode and the solution has exponential integrability.
In view of Theorem \ref{th1}, to show the Sobolev regularity of the strong solution, we will pay our attention
on the $p$-th moment estimates of $X_t(x)-X_t(y)$. This is the place where assumptions (\ref{c4})-(\ref{cf}) are needed.
As in \cite{Zh3,Zh1}, the estimates of Krylov's type will play an important role throughout this paper.
However, since we are assuming only some local integrability conditions and the coefficients may have exponential growth rate at infinity,
some new probabilistic estimates are established (see Lemma \ref{m} and Lemmas \ref{Le33}, \ref{Le34} below).

\vspace{2mm}

To illustrate Theorem \ref{main1}, we present below two examples.
\bx
Consider the following one-dimensional SDE:
$$
\dif X_t=[(1-X^5_t)1_{X_t<0}-(1+X^5_t)1_{X_t\geq0}]\dif t+(1+|X_t|^2)^\beta\dif W_t,
$$
where $\beta\in[0,1)$.
In this case, $\sigma(x)=(1+x^2)^\beta$ and $b(x)=(1-x^5)1_{x<0}-(1+x^5)1_{x\geq 0}$  are both of super-linear growth,
and the drift $b$ has a jump at $0$. Moreover, for any $\kappa>0$, by Young's inequality, it is easy to see that
$$
\<x,b(x)\>+\kappa(1+x^2)|\sigma(x)|^2=-x^6-|x|+\kappa(1+x^2)^{1+2\beta}\leq C_{\kappa,\beta}
$$
and $\<x-y,b(x)-b(y)\>\leq 0$,
$$
 \|\sigma(x)-\sigma(y)\|^2\leq C_{\beta} |x-y|^2\Big(1+|x|^{2(2\beta-1)\vee 0}+|y|^{2(2\beta-1)\vee 0}\Big).
$$
Thus, {\bf (H1)}, {\bf (H2)} and \eqref{Sigma} hold.
\ex
\bx
Suppose that for any $\kappa>0$ and $T>0$, there is a convex function $F_\kappa(x)$ such that
$$
\sup_{|\xi|=1}\<\xi, \nabla_\xi b(t,x)\>+\kappa\|\nabla\sigma(t,x)\|^2\leq F_\kappa(x),
$$
where $\nabla_\xi f:=\<\nabla f,\xi\>$ for a $C^1$-function $f:\mR^d\to\mR$. Under this assumption, (\ref{c4}) holds. In fact, by the mean-value formula, we have
\begin{align*}
&(\<x-y,b(t,x)-b(t,y)\>+\kappa\|\sigma(t,x)-\sigma(t,y)\|^2)/|x-y|^2\\
&\leq \int^1_0\Big[\sup_{|\xi|=1}\<\xi, \nabla_\xi b(t,\theta x+(1-\theta) y)\>+\kappa\|\nabla\sigma(t,\theta x+(1-\theta) y)\|^2\Big]\dif\theta\\
&\leq\int^1_0F_\kappa(\theta x+(1-\theta) y)\dif\theta\leq \frac{F_\kappa(x)+F_\kappa(y)}{2},
\end{align*}
where the last step is due to the convexity of $F_\kappa$.
Compared with \cite{Ch-Li}, our assumptions {\bf (H1)} and {\bf (H2)} are significantly weaker.
\ex


In Theorem \ref{main1}, the drift $b$ is locally bounded.
Our next result allows the drift $b$ to be locally singular and of linear growth at infinity. To this aim, we make the following assumptions:
\vspace{1mm}
\begin{enumerate}[{\bf(H1$'$)}]
\item {\bf (Local Sobolev integrability)} $\sigma$ is uniformly continuous in $x$ and locally uniformly with respect to $t\in\mR_+$,
and for some $q>2d+2$,
$$
b\in L^{q}_{loc}(\mR_+\times\mR^d),\,\, \nabla\sigma\in L^{q}_{loc}(\mR_+\times\mR^d),
$$
and for any $T>0$, there is a constant $K\geq 1$ such that  for all $(t,x)\in[0,T]\times\mR^d$,
$$
K^{-1} |\xi|^2\leq |\sigma(t,x)\xi|^2\leq K|\xi|^2,\quad\forall\xi\in\mR^m.
$$
\item {\bf (Lipschitz continuity outside a ball)} For any $T>0$, there exist $R_0\geq 1$, $\alpha'\in[0,1)$ and constant $C>0$ such that for all $t\in[0,T]$,
$$
|b(t,x)|\leq C(1+|x|),\ \ |x|\geq R_0,
$$
and for all $t\in[0,T]$ and $|x|, |y|\geq R_0$,
\begin{align}
\begin{split}
&\qquad\|\sigma(t,x)-\sigma(t,y)\|\leq C|x-y|, \\
&|b(t,x)-b(t,y)|\leq C|x-y|\Big(|x|^{2\alpha'}+|y|^{2\alpha'}\Big). 
\end{split}
\label{c041}
\end{align}
\end{enumerate}

It should be noticed that conditions in ${\bf{(H2')}}$ are assumed to hold only outside a large ball, while $b$ can be singular in the ball.
We have
\bt\label{main3}
Under ${\bf{(H1')}}$ and ${\bf{(H2')}}$, there exists a unique global strong solution $X_t(x)$ to SDE (\ref{1}) so that
$(t,x)\mapsto X_t(x)$ is continuous. Moreover, the conclusions {\bf (A)} with \eqref{so1}, {\bf (C)} and {\bf (D)} in Theorem \ref{main1} still hold,
and the $\alpha$ in \eqref{so1} is taken to be $1$.
\et
We organize this paper as follows: In Section 2, we make some preparations, and give the proof of Theorem \ref{th1}
and a criterion on the existence of exponential moments of a Markov process.
In Section 3, we provide some estimates on the solution to equation (\ref{1}) and give the proof of Theorem \ref{main1}.
Finally, the proof of Theorem \ref{main3} is given in Section 4 by using Zvonkin's transformation and Theorem \ref{main1}.

\vspace{3mm}

Throughout this paper, we use the following convention: $C$ with or without subscripts will denote a positive constant, whose value may change in different places, and
whose dependence on the parameters can be traced from the calculations.

\section{Preliminaries}

Let $U$ be an open domain in $\mR^d$. For $p\in[1,\infty]$, let $\mW^{1,p}(U)$ be the classical first order Sobolev space:
$$
\mW^{1,p}(U):=\Big\{f\in L^1_{loc}(U): \|f\|_{1,p}:=\|f\|_{p}+\|\nabla f\|_{p}<+\infty\Big\},
$$
where $\|\cdot\|_p$ is the usual $L^p(U)$-norm and $\nabla$ denotes the gradient in the distributional sense.
When $U$ is a bounded $C^1$-domain, it was proved in \cite{Ha} that a function $f\in\mW^{1,p}(U)$ if and
only if $f\in L^p(U)$ and there exists a nonnegative function $g\in L^p(U)$ such that for Lebesgue-almost all $x,y\in U$,
$$
|f(x)-f(y)|\leq |x-y|\big(g(x)+g(y)\big).
$$

Let us now extend the above characterization to the case of random fields.
For $p,q,r\in[1,\infty]$ and $T>0$, let $L^r(T):=L^r([0,T])$ and define
\begin{align}
\begin{split}
&\mW^{1,q}\big(U;L^p(\Omega;L^r(T))\big)\\
&\quad:=\Big\{f\in L^1_{loc}(U\times[0,T]; L^1(\Omega)):f,\,\nabla f\in L^{q}\big(U;L^p(\Omega; L^r(T))\big)\Big\},
\end{split} \label{w}
\end{align}
and
$$
\|f\|_{\mW^{1,q}(U;L^p(\Omega;L^r(T)))}:=\|f\|_{L^{q}(U;L^p(\Omega;L^r(T)))}+\|\nabla f\|_{L^{q}(U;L^p(\Omega;L^r(T)))}.
$$
Notice that by Fubini's theorem,
$$
L^p\big(U;L^p([0,T]\times\Omega)\big)=L^{p}\big([0,T]\times\Omega;L^p(U)\big),
$$
and hence,
\begin{align}
\mW^{1,p}\big(U;L^p([0,T]\times\Omega)\big)=L^p\big([0,T]\times\Omega;\mW^{1,p}(U)\big).\label{RG5}
\end{align}
In what follows, we write $B_r:=\{x\in\mR^d: |x|<r\}$.
\bl\label{Le21}
Let $p\in(1,\infty)$, $q,r\in(1,\infty]$ and $f\in L^{q}\big(U;L^p(\Omega; L^r(T))\big)$. Assume that there exists a nonnegative measurable function $g\in L^q(U)$ such that
for Lebesgue-almost all $x,y\in U$,
\begin{align}
\|f(x,\cdot)-f(y,\cdot)\|_{L^p(\Omega; L^r(T))}\leq |x-y|\big(g(x)+g(y)\big), \label{Mo}
\end{align}
then $f\in\mW^{1,q}\big(U;L^p(\Omega; L^r(T))\big)$, and for Lebesgue-almost all $x\in U$,
\begin{align}\label{UY8}
\|\p_if(x,\cdot)\|_{L^p(\Omega; L^r(T))}\leq 2g(x),\ \ i=1,\cdots,d.
\end{align}
\el

\begin{proof}
Below, we always extend a function $f$ defined on $U$ to $\mR^d$ by setting $f(x,\cdot)\equiv0$ for $x\notin U$.
Let $\varrho:\mR^d\to[0,1]$ be a smooth function with support in $B_1$ and
$\int\varrho\dif x=1$. For $n\in\mN$, define a family of mollifiers $\varrho_n(x)$ as follows:
\begin{align}
\varrho_n(x):=n^{d}\varrho(nx).\label{Rho}
\end{align}
Define the mollifying approximations of $f$ and $g$ by
\begin{align}
f_n(x,t,\omega):=f(\cdot,t,\omega)*\varrho_n(x),\ \ g_n(x):=g*\varrho_n(x).\label{fn}
\end{align}
For $\eps\in(0,1]$, set
$$
U_\eps:=\{x\in U: \dis(x,\p U)>\eps\},
$$
where $\dis(x,\p U)$ denotes the distance between $x$ and the boundary $\p U$.
By (\ref{Mo}), it is easy to see that for any $x,y\in U_\eps$ and $n>2/\eps$,
\begin{align}
\begin{split}
\|f_n(x)-f_n(y)\|_{L^p(\Omega; L^r(T))}&\leq\int_{\mR^d}\|f(x-z)-f(y-z)\|_{L^p(\Omega; L^r(T))}\varrho_n(z)\dif z\\
&\leq|x-y|\int_{\mR^d}(g(x-z)+g(y-z))\varrho_n(z)\dif z\\
&=|x-y|(g_n(x)+g_n(y)).
\end{split}\label{Mo1}
\end{align}
Let $\{e_i,i=1,\cdots,d\}$ be the canonical basis of $\mR^d$. For all $x\in U_{\eps}$ and $n>2/\eps$, by Fatou's lemma and (\ref{Mo1}),
we have
\begin{align}
\begin{split}
\|\p_if_n(x)\|_{L^p(\Omega; L^r(T))}&=\left\|\lim_{\delta\to 0}\frac{|f_n(x+\delta e_i)-f_n(x)|}{\delta}\right\|_{L^p(\Omega; L^r(T))}\\
&\leq\varliminf_{\delta\to 0}\frac{\|f_n(x+\delta e_i)-f_n(x)\|_{L^p(\Omega; L^r(T))}}{\delta}\\
&\leq\varliminf_{\delta\to 0}(g_n(x+\delta e_i)+g_n(x))=2g_n(x).
\end{split}\label{RG7}
\end{align}
Integrating both sides on $U_{\eps}$ , we obtain
\begin{align}
\int_{U_{\eps}}\|\p_if_n(x)\|^q_{L^p(\Omega; L^r(T))}\dif x\leq 2^q\int_{U_{\eps}}g_n(x)^q\dif x\leq 2^q\|g\|^q_{L^q(U)}.\label{RG6}
\end{align}
In particular, if we let $\gamma=p\wedge q\wedge r$ and $U^R_\eps:=U_{\eps}\cap B_R$ for $R>0$, then by (\ref{RG5}), we have for any $R\in\mN$,
$$
\sup_n\|f_n\|_{L^\gamma([0,T]\times\Omega; \mW^{1,\gamma}(U^R_{\eps}))}=
\sup_n\|f_n\|_{\mW^{1,\gamma}(U^R_{\eps}; L^\gamma([0,T]\times\Omega))}<\infty.
$$
Since $\gamma\in(1,\infty)$ and $L^\gamma([0,T]\times\Omega; \mW^{1,\gamma}(U^R_\eps))$ is weakly compact
and $f_n\to f$ in $L^\gamma([0,T]\times\Omega; L^\gamma(U))$,
we have
\begin{align}
f\in L^\gamma([0,T]\times\Omega; \mW^{1,\gamma}(U^R_{\eps})).\label{RG8}
\end{align}
By the arbitrariness of $\eps$ and $R$, one sees that for $\dif t\times\mP(\dif\omega)$-almost all $(t,\omega)$, $x\mapsto f(x,t,\omega)$ is weakly differentiable in $U$, and
for all $x\in U_{\eps}$ and $n>2/\eps$,
$$
\p_i f_n(x,t,\omega)=\p_i f(\cdot,t,\omega)*\varrho_n(x),
$$
which, by (\ref{RG8}) and the property of convolutions, then implies that
$$
\lim_{n\to\infty}\|\p_i f_n-\p_i f\|_{L^\gamma(U^R_{\eps}\times[0,T]\times\Omega)}=0.
$$
Thus, for some subsequence $n_k$ and $\dif x\times\dif t\times\mP(\dif\omega)$-almost all $(x,t,\omega)\in U^R_\eps\times[0,T]\times\Omega$,
$$
\p_i f_{n_k}(x,t,\omega)\to\p_if(x,t,\omega).
$$
Now, by (\ref{RG7}) and Fatou's lemma, we obtain that for Lebesgue-almost all $x\in U^R_\eps$,
$$
\|\p_if(x)\|_{L^p(\Omega; L^r(T))}\leq\varliminf_{k\to\infty}\|\p_if_{n_k}(x)\|_{L^p(\Omega; L^r(T))}\leq 2\varliminf_{k\to\infty}g_{n_k}(x)=2g(x).
$$
The proof is complete by the arbitrariness of $\eps$ and $R$.
\end{proof}

We also have the following converse result.
\bl\label{Le22}
Let $f\in\mW^{1,q}_{loc}\big(\mR^d; L^p(\Omega; L^r(T))\big)$ for some $p,q,r\in(1,\infty]$. For any $R>0$, there exists a measurable function
$g_R\in L^q_{loc}(\mR^d)$ such that for Lebesgue-almost all $x,y\in\mR^d$ with $|x-y|<R$,
\begin{align}
\|f(x,\cdot)-f(y,\cdot)\|_{L^p(\Omega; L^r(T))}\leq |x-y|\big(g_R(x)+g_R(y)\big). \label{Mo2}
\end{align}
Moreover, if $f\in\mW^{1,q}\big(\mR^d; L^p(\Omega; L^r(T))\big)$, then $R$ can be $\infty$ and $g_\infty\in L^q(\mR^d)$.
\el
\begin{proof}
Let $f_n$ be the mollifying approximation of $f$ as in (\ref{fn}). By \cite[Lemma 3.5]{Zh0}, we have
\begin{align*}
|f_n(x,t,\omega)-f_n(y,t,\omega)|&\leq 2^d\int^{|x-y|}_0\!\!\!\fint_{B_s}|\nabla f_n(x+z,t,\omega)|\dif z\dif s\\
&+2^d\int^{|x-y|}_0\!\!\!\fint_{B_s}|\nabla f_n(y+z,t,\omega)|\dif z\dif s.
\end{align*}
Hence, for all $x,y\in\mR^d$ with $|x-y|<R$,
\begin{align}
\begin{split}
\|f_n(x)-f_n(y)\|_{L^p(\Omega; L^r(T))}&\leq 2^d\int^{|x-y|}_0\!\!\!\fint_{B_s}\|\nabla f_n(x+z)\|_{L^p(\Omega; L^r(T))}\dif z\dif s\\
&+2^d\int^{|x-y|}_0\!\!\!\fint_{B_s}\|\nabla f_n(y+z)\|_{L^p(\Omega; L^r(T))}\dif z\dif s\\
&\leq 2^d|x-y|(g_R^n(x)+g^n_R(y)),
\end{split}\label{Lim}
\end{align}
where
$$
g^n_R(x):=\cM_R\|\nabla f_n\|_{L^p(\Omega; L^r(T))}(x):=\sup_{s\in(0,R)}\fint_{B_s}\|\nabla f_n(x+z)\|_{L^p(\Omega; L^r(T))}\dif z
$$
is the local maximal function of $x\mapsto\|\nabla f_n(x)\|_{L^p(\Omega; L^r(T))}$. Notice that
\begin{align*}
\cM_R\|\nabla f_n\|_{L^p(\Omega; L^r(T))}(x)&\leq\sup_{s\in(0,R)}\fint_{B_s}\|\nabla f\|_{L^p(\Omega; L^r(T))}*\varrho_n(x+z)\dif z\\
&\leq\cM_R\|\nabla f\|_{L^p(\Omega; L^r(T))}*\varrho_n(x),
\end{align*}
and $\cM_R\|\nabla f\|_{L^p(\Omega; L^r(T))}\in L^q_{loc}(\mR^d)$ by the property of maximal functions (cf. \cite{St}).
By taking limits for both sides of
\begin{align*}
&\|f_n(x)-f_n(y)\|_{L^p(\Omega; L^r(T))}\\
&\quad\leq 2^d|x-y|\Big(\cM_R\|\nabla f\|_{L^p(\Omega; L^r(T))}*\varrho_n(x)+\cM_R\|\nabla f\|_{L^p(\Omega; L^r(T))}*\varrho_n(y)\Big),
\end{align*}
we obtain that for Lebesgue-almost all $x,y\in\mR^d$ with $|x-y|<R$,
\begin{align}\label{UY5}
\begin{split}
&\|f(x)-f(y)\|_{L^p(\Omega; L^r(T))}\\
&\quad\leq 2^d|x-y|\Big(\cM_R\|\nabla f\|_{L^p(\Omega; L^r(T))}(x)+\cM_R\|\nabla f\|_{L^p(\Omega; L^r(T))}(y)\Big).
\end{split}
\end{align}
The proof is complete.
\end{proof}

Combining Lemma \ref{Le21} with Lemma \ref{Le22}, we can give

\begin{proof}[Proof of Theorem \ref{th1}]
The sufficiency follows by Lemma \ref{Le21}. Let 
$$
f\in\mW^{1,q}\big(U;L^p(\Omega; L^r(T))\big).
$$
Since $U$ is a bounded $C^1$-domain, there exists an extension operator  (cf. \cite[p.151, Theorem 5.22]{Ad} or \cite[Chapter VI]{St})
$$
\mT:\mW^{1,q}\big(U;L^p(\Omega; L^r(T))\big)\to \mW^{1,q}\big(\mR^d;L^p(\Omega; L^r(T))\big)
$$
such that $\mT f=f$ restricted on $U$ and
$$
\|\mT f\|_{\mW^{1,q}(\mR^d;L^p(\Omega; L^r(T)))}\leq C\|f\|_{\mW^{1,q}(U;L^p(\Omega; L^r(T)))}.
$$
Thus, (\ref{Mo23}) follows by (\ref{Mo2}).
\end{proof}


We also need the following local version of Khasminskii's estimate (see \cite[Lemma 2.1]{Sz}).

\bl\label{m}
Let $(\Omega,\sF,(\mP_x)_{x\in\mR^d}; (X_t)_{t\geq 0})$ be a family of $\mR^d$-valued time-homogenous
Markov process. Let $f$ be a nonnegative measurable function over $\mR^d$. For given $T,R>0$, if
\begin{align}
\sup_{|x|\leq R}\mE_{x}\Bigg(\int_0^T\!\!f(X_t)1_{|X_t|\leq R}\dif t\Bigg)=:c<1,  \label{y1}
\end{align}
where $\mE_x$ denotes the expectation with respect to $\mP_x$, then for all $x\in\mR^d$,
\begin{align}
\mE_{x}\exp\Bigg\{\int_0^T\!f(X_t)1_{|X_t|\leq R}\dif t\Bigg\}\leq 1+\frac{1}{1-c}\mE_x\left(\int_0^T\!f(X_t)1_{|X_t|\leq R}\dif t\right).   \label{y2}
\end{align}
\el
\begin{proof}
Set $f_R(x):=f(x)1_{|x|\leq R}$. By Taylor's expansion, we can write
$$
\mE_x\exp\Bigg\{\int_0^Tf_R(X_t)\dif t\Bigg\}=\sum_{n=0}^{\infty}\frac{1}{n!}\mE_x\Bigg(\int_0^Tf_R(X_t)\dif t\Bigg)^n.
$$
For $n\in\mN$, noticing that
$$
\left(\int^T_0g(t)\dif t\right)^n=n!\int\!\!\!...\!\!\!\int_{\Delta_T^n}g(t_1)\cdots g(t_n)\dif t_1\cdots\dif t_n,
$$
where
$$
\Delta_T^n:=\Big\{(t_1,\cdots, t_n): 0\leq t_1\leq t_2\leq\cdots\leq t_n\leq T\Big\},
$$
we further have
\begin{align*}
&\mE_x\exp\Bigg\{\int_0^Tf_R(X_t)\dif t\Bigg\}
=1+\sum_{n=1}^{\infty}\mE_x\Bigg(\int\!\!\!...\!\!\!\int_{\Delta_T^n}f_R(X_{t_{1}})\cdots f_R(X_{t_{n}})\dif t_1\cdots\dif t_n\Bigg)\\
&=1+\sum_{n=1}^{\infty}\mE_x\Bigg(\int\!\!\!...\!\!\!\int_{\Delta_T^{n-1}}f_R(X_{t_{1}})\cdots f_R(X_{t_{n-1}})
\mE_{X_{t_{n-1}}}\int_0^{T-t_{n-1}}f_R(X_{t_{n}})\dif t_n\dif t_1\cdots\dif t_{n-1}\Bigg)\\
&\!\!\!\stackrel{(\ref{y1})}{\leq}1+\sum_{n=1}^{\infty}c\mE_x\Bigg(\int\!\!\!...\!\!\!\int_{\Delta_T^{n-1}}f_R(X_{t_{1}})\cdots f_R(X_{t_{n-1}})\dif t_1\cdots\dif t_{n-1}\Bigg)
\leq\cdots\\
&\leq1+\sum_{n=1}^{\infty}c^{n-1}\mE_x\Bigg(\int_0^Tf_R(X_{t_1})\dif t_1\Bigg)=1+\frac{1}{1-c}\mE_x\Bigg(\int_0^Tf_R(X_{t})\dif t\Bigg),
\end{align*}
where the second equality is due to the Markov property of $X_t$.
\end{proof}
Finally, we recall the following Krylov estimate about the distributions of continuous semimartingales (cf. \cite{Kry} or \cite[Lemma 3.1]{Gy-Ma}).
\bl\label{kry}
Let $m=m_t$ be a continuous $\mR^d$-valued local martingale, and $V=V_t$ a continuous $\mR^d$-valued process with finite variation on finite time intervals.
Suppose that
$$
m(0) = V(0) = 0,\ \ \dif \<m\>_t\ll\dif t,
$$
and set
$$
a(t):=\tfrac{\dif\<m\>_t}{2\dif t},\ \ X(t):=m(t)+ V(t).
$$
For any $\lambda> 0$, stopping time $\tau$ and nonnegative Borel function
$f: \mR_+\times\mR^d\rightarrow\mR_+$, we have
\begin{align}
\begin{split}
&\mE\int_0^{\tau}\e^{-\lambda t}\big(\det a(t)\big)^{\frac{1}{d+1}}f(t,X_t)\dif t\\
&\quad\leq N_{d,\lambda}(\mV^2+\mA)^{\frac{d}{2(d+1)}}\Bigg(\int_0^{\infty}\!\!\!\int_{\mR^d}f^{d+1}(t,x)\dif x\dif t\Bigg)^{\frac{1}{d+1}},   
\end{split}
\label{k1}
\end{align}
where
\begin{align}
\mV:=\mE\int_0^{\tau}\e^{-\lambda t}|\dif V(t)|,\, \mA:=\mE\int_0^{\tau}\e^{-\lambda t}\t\,a(t)\dif t,\label{RG1}
\end{align}
and $N_{d,\lambda}$ is a constant depending only on $d$ and $\lambda$.
\el

\section{Proof of Theorem \ref{main1}}
Below we write
$$
\sL^{\sigma,b}_s f(x):=\tfrac{1}{2}\sum_{ijk}\sigma_{ik}(s,x)\sigma_{jk}(s,x)\p_i\p_j f(x)+\sum_ib_i(s,x)\p_i f(x).
$$
Under ${\bf{(H1)}}$, it has been proven in \cite{Zh1, Zh4} that SDE (\ref{1}) admits a unique local strong solution.
The following lemma gives the non-explosion and the exponential integrability of $X_t(x)$ under \eqref{c1}.

\bl\label{Le31}
Let $X_t(x)$ be the unique local solution of (\ref{1}) with starting point $x$. Under ${\bf{(H1)}}$ and (\ref{c1}), there is a unique global solution $X_t(x)$ to SDE \eqref{1}.
Moreover, let $\alpha\in[0,1]$ and $\kappa\mapsto C_{\kappa}\in\mR$ be as in \eqref{c1}.
\begin{itemize}
\item($\alpha>0$) For any $\lambda\geq 2\alpha C_{\alpha+1}$ and for all $t\geq0$ and $x\in\mR^d$, we have
\begin{align}
\mE\exp\Big\{\e^{-\lambda t}\big(1+|X_t(x)|^2\big)^{\alpha}\Big\}\leq \exp\Big\{(1+|x|^2)^{\alpha}\Big\}.   \label{ex}
\end{align}
\item ($\alpha=0$) For any $p\geq 1$ and $\lambda\geq C_p$ and for all $t\geq 0$ and $x\in\mR^d$, we have
\begin{align}
\mE(1+|X_t(x)|^2)^{p}\leq \e^{\lambda t}(1+|x|^2)^p.  \label{ex1}
\end{align}
\end{itemize}
\el
\begin{proof}
We only consider the case of $\alpha>0$. For $\alpha=0$, it is similar. For $R>0$, define
\begin{align}\label{TAU}
\tau_R:=\inf\Big\{t\geq 0: |X_t(x)|\geq R\Big\},
\end{align}
and for $\lambda\in\mR$,
\begin{align}
f(t,x):=\exp\Big\{\e^{-\lambda t}(1+|x|^2)^{\alpha}\Big\}.\label{EB2}
\end{align}
By It\^o's formula, we have
\begin{align*}
\mE f(t\wedge\tau_R,X_{t\wedge\tau_R})=f(0,x)+\mE\left(\int^{t\wedge\tau_R}_0(\p_s f+\sL^{\sigma,b}_s f)(s,X_s)\dif s\right).
\end{align*}
Notice that
\begin{align}
(\p_s f+\sL^{\sigma,b}_s f)(s,x)&=\alpha(1+|x|^2)^\alpha\e^{-\lambda s}f(s,x)
\Big(-\tfrac{\lambda}{\alpha}+\tfrac{2\<b,x\>+\|\sigma\|^2}{1+|x|^2}\no\\
&+2\sum_{ijk}\sigma_{ik}\sigma_{jk}\left[\alpha(1+|x|^2)^{\alpha-2}\e^{-\lambda s}
+\tfrac{\alpha-1}{(1+|x|^2)^2}\right]x_ix_j\Big)\no\\
&\leq\alpha(1+|x|^2)^\alpha\e^{-\lambda s}f(s,x)\left(-\tfrac{\lambda}{\alpha}+\tfrac{2\<b,x\>+2(\alpha+1)(1+|x|^2)^\alpha\|\sigma\|^2}{1+|x|^2}\right)\label{EB1}\\
&\stackrel{(\ref{c1})}{\leq}\alpha(1+|x|^2)^\alpha\e^{-\lambda s}f(s,x)\left(-\tfrac{\lambda}{\alpha}+2C_{\alpha+1}\right).\no
\end{align}
Hence, if $\lambda\geq 2\alpha C_{\alpha+1}$, then
$$
\mE f(t\wedge\tau_R,X_{t\wedge\tau_R})\leq f(0,x).
$$
By letting $R\to\infty$, one sees that $\tau_\infty=\infty$, i.e., no explosion, and (\ref{ex}) holds.
\end{proof}

The following global Krylov estimate is an easy consequence of Lemma \ref{kry} and Lemma \ref{Le31}.
\bl\label{Le32}
Under {\bf (H1)}, (\ref{c1}) and (\ref{c3}), for any $q>d+1$ and $T>0$, there exist constants $C, \gamma>0$
such that for all nonnegative $f\in L^q([0,T]\times\mR^d)$ and $x\in\mR^d$,
$$
\mE\left(\int_0^Tf(t,X_t(x))\dif t\right)\leq C\Big(1_{\alpha>0}\e^{(1+|x|^2)^\alpha}+1_{\alpha=0}(1+|x|^2)^{\gamma}\Big)
\left(\int_0^T\!\!\!\int_{\mR^d}f^{q}(t,y)\dif y\dif t\right)^{\frac{1}{q}}.
$$
\el
\begin{proof}
In Lemma \ref{kry}, let us take
$$
m(t):=\int^t_0\sigma(s,X_s)\dif W_s,\ \ V(t):=\int^t_0b(s,X_s)\dif s,
$$
so that
$$
a(t)=\frac{\dif \<m\>_t}{2\dif t}=\tfrac{1}{2}(\sigma\sigma^*)(t,X_t).
$$
By H\"older's inequality and Lemma \ref{kry}, we have
\begin{align}
\begin{split}
&\mE\left(\int_0^Tf(t,X_t(x))\dif t\right)=\mE\int_0^{T}f(t,X_t(x))\big(\det a(t)\big)^{\frac{1}{q}}\big(\det a(t)\big)^{-\frac{1}{q}}\dif t\\
&\leq\Bigg(\mE\int_0^{T}f^{\frac{q}{d+1}}(t,X_t(x))\big(\det a(t)\big)^{\frac{1}{d+1}}\dif t\Bigg)^{\frac{d+1}{q}}
\Bigg(\mE\int_0^{T}\!\!\big(\det a(t)\big)^{-\frac{1}{q-d-1}}\dif t\Bigg)^{\frac{q-d-1}{q}}\\
&\leq C_{T,d}(\mV^2+\mA)^{\frac{d}{2q}}\Bigg(\int_0^T\!\!\!\int_{\mR^d}f^q(t,x)\dif x\dif t\Bigg)^{\frac{1}{q}}
\Bigg(\mE\int_0^{T}\!\!\big(\det a(t)\big)^{-\frac{1}{q-d-1}}\dif t\Bigg)^{\frac{q-d-1}{q}},
\end{split}\label{TR3}
\end{align}
where $\mA$ and $\mV$ are defined by (\ref{RG1}).
By (\ref{c3}), (\ref{ex}) and Young's inequality, we have
\begin{align}
\begin{split}
\mA+\mV^2&\leq C\mE\int_0^T(|b(t,X_t)|^2+\|\sigma(t,X_t)\|^2)\dif t\\
&\leq C\mE\int_0^T\Big(1_{\alpha>0}\exp\{2C_2(1+|X_t|^2)^{\alpha'}\}+1_{\alpha=0}C^2_2(1+|X_t|^2)^{2\gamma_2}\Big)\dif t\\
&\leq C\int_0^T\Big(1_{\alpha>0}\mE\exp\{\e^{-\lambda t}(1+|X_t|^2)^{\alpha}\}+1_{\alpha=0}\mE(1+|X_t|^2)^{2\gamma_2}\Big)\dif t\\
&\leq C\Big(1_{\alpha>0}\e^{(1+|x|^2)^\alpha}+1_{\alpha=0}(1+|x|^2)^{2\gamma_2}\Big).
\end{split}\label{TR1}
\end{align}
Similarly, by (\ref{c2}), we have for some $\gamma_3>0$,
\begin{align}
\mE\int_0^{T}\!\!\big(\det a(t)\big)^{-\frac{1}{q-d-1}}\dif t\leq C\Big(1_{\alpha>0}\e^{(1+|x|^2)^\alpha}+1_{\alpha=0}(1+|x|^2)^{\gamma_3}\Big).\label{TR2}
\end{align}
Substituting (\ref{TR1}) and (\ref{TR2}) into (\ref{TR3}), we obtain the desired estimate.
\end{proof}

Taking into account Lemma \ref{m}, we can prove the following global-exponential moment estimate of Krylov's type,
which will play a crucial role in the proof of Theorem \ref{main1}. Since $f$ is allowed to be singular in a ball and of linear growth at infinity,
we need to separately consider the interior and exterior parts of a ball by using Lemmas \ref{m}, \ref{Le31} and \ref{Le32}.
\bl\label{Le33}
For given $q>d+1$, let $f\in L^q_{loc}(\mR_+\times\mR^d)$ be a nonnegative measurable function.
Let $\alpha$ be as  in (\ref{c1}). Suppose that for some $R_0, C_0>0$ and $\alpha'\in[0,\alpha)$,
\begin{align}
f(t,x)1_{\{|x|>R_0\}}\leq C_0\Big[1_{\alpha>0}(1+|x|^2)^{\alpha'}+1_{\alpha=0}\log(1+|x|^2)\Big].\label{RG2}
\end{align}
Under {\bf (H1)}, (\ref{c1}) and (\ref{c3}), for any  $T>0$, there are $C,\gamma>0$ such that for all $x\in\mR^d$,
\begin{align}
\mE\exp\left\{\int^T_0f(t,X_t(x))\dif t\right\}\leq C\Big(1_{\alpha>0}\e^{(1+|x|^2)^\alpha}+1_{\alpha=0}(1+|x|^2)^{\gamma}\Big).\label{ER1}
\end{align}
\el
\begin{proof}
Set
$$
f_{R_0}(t,x):=f(t,x)1_{|x|\leq R_0},\ \ \bar f_{R_0}(t,x):=f(t,x)1_{|x|>R_0}.
$$
By H\"older's inequality, we have
\begin{align}
\begin{split}
&\left(\mE\exp\left\{\int^T_0f(t,X_t(x))\dif t\right\}\right)^2\leq\mE\exp\left\{2\int^T_0\bar f_{R_0}(t, X_t(x))\dif t\right\}\\
&\qquad\times\mE\exp\left\{2\int^T_0f_{R_0}(t, X_t(x))\dif t\right\}=:I_1(T,x)\times I_2(T,x).
\end{split}\label{UY2}
\end{align}
For $I_1(T,x)$, by (\ref{RG2}), Jensen's inequality and Lemma \ref{Le31}, we have
\begin{align}
I_1(T,x)&\leq \mE\exp\left\{2C_0\int^T_0\Big[1_{\alpha>0}(1+|X_t(x)|^2)^{\alpha'}+1_{\alpha=0}\log(1+|X_t(x)|^2)\Big]\dif t\right\}\no\\
&\leq \frac{1}{T}\int^T_0 \mE\exp\left\{2C_0 T\Big[1_{\alpha>0}(1+|X_t(x)|^2)^{\alpha'}+1_{\alpha=0}\log(1+|X_t(x)|^2)\Big]\right\}\dif t\no\\
&\leq \frac{C}{T}\int^T_0 \left(1_{\alpha>0}\mE\exp\left\{\e^{-\lambda t}(1+|X_t(x)|^2)^{\alpha}\right\}+1_{\alpha=0}\mE(1+|X_t(x)|^2)^{2C_0 T}\right)\dif t\label{UY1}\\
&\leq C\left(1_{\alpha>0}\exp\Big\{(1+|x|^2)^\alpha\Big\}+1_{\alpha=0}(1+|x|^2)^{2C_0T}\right),\no
\end{align}
where $\lambda$ is the same as in \eqref{ex}, and the third inequality is due to Young's inequality.
For $I_2(T,x)$, for any $\eps\in(0,1)$ and $\theta>1$, by Young's inequality we have
\begin{align}\label{UY3}
I_2(T,x)\leq \e^{C_\eps}\mE\exp\left\{\eps\int^T_0f_{R_0}(t, X_t(x))^\theta\dif t\right\}.
\end{align}
Let us choose $\theta>1$ so that $\frac{q}{\theta}>d+1$. Then by Lemma \ref{Le32}, we have
$$
\eps\sup_{|x|\leq R_0}\mE\left(\int^T_0f_{R_0}(t, X_t(x))^\theta\dif t\right)
\leq \eps C_{R_0}\left(\int^T_0\!\!\!\int_{|y|<R_0}f(t,y)^{q}\dif y\dif t\right)^{\frac{\theta}{q}}=:c_\eps.
$$
Since $(t,X_t(x))$ is a time-homogenous Markov process in $\mR_+\times\mR^d$,
if we choose $\eps$ being small enough so that $c_\eps<1$, then by (\ref{y2}), we obtain
$$
\mE\exp\left\{\eps\int^T_0f_{R_0}(t, X_t(x))^\theta\dif t\right\}
\leq 1+\frac{\eps}{1-c_\eps}\mE\left(\int^T_0f_{R_0}(t, X_t(x))^\theta\dif t\right),
$$
which, together with \eqref{UY2}, \eqref{UY1}, \eqref{UY3} and Lemma \ref{Le32}, yields the desired estimate.
\end{proof}

We also need the following local-exponential moment estimate of Krylov's type,
which is a consequence of Lemma \ref{m} and \cite[Theorem 2.1]{Zh1} (see also \cite[Theorem 4.1]{Zh4}).
In particular, the integrability index $q$ in the following lemma can be smaller than the one in Lemma \ref{Le33}.
\bl\label{Le34}
For $R\geq 1$, let $\tau_R$ be defined by \eqref{TAU}. Under {\bf (H1)} and (\ref{c3}), for any $q>\frac{d}{2}+1$ and $T>0$, there exists a constant $C_R>0$
such that for all $f\in L^q_{loc}(\mR^{d+1})$ and $|x|<R$,
\begin{align}
\mE\exp\left\{\int_0^{T\wedge\tau_R}f(t,X_t(x))\dif t\right\}<\infty.\label{k22}
\end{align}
\el
\begin{proof}
Let $\chi_R$ be a smooth cutoff function with $\chi_R(x)=1$ for $|x|\leq R$ and $\chi_{R}=0$ for $|x|\geq R+1$, and set
$$
b_R(t,x)=b(t,x)\chi_R(x),\ \ \sigma_R(t,x):=\sigma(t,\chi_R(x)x).
$$
By {\bf (H1)} and (\ref{c3}), it is easy to see that for some $C_R>0$ and $\alpha\in(0,1)$,
$$
|b_R(t,x)|\leq C_R,\ \ C_R^{-1}|\xi|\leq |\sigma_R(t,x)\xi|\leq C_R|\xi|,\ \ |\sigma_R(t,x)-\sigma_R(t,y)|\leq C_R|x-y|^\alpha.
$$
Let $X^R_t(x)$ solve SDE \eqref{1} with $(b_R,\sigma_R)$ in place of $(b,\sigma)$. By the local uniqueness, one has
$$
X^R_t(x)=X_t(x),\ \ t<\tau_R.
$$
Hence, letting $\theta>1$ so that $\frac{q}{\theta}>\frac{d}{2}+1$, by \cite [Theorem 2.1]{Zh1}, we have for all $|x|<R$,
\begin{align*}
&\mE\left(\int_0^{T\wedge\tau_R}|f(t,X_t(x))|^\theta\dif t\right)
=\mE\left(\int_0^{T}|f(t,X^R_t(x))|^\theta 1_{t<\tau_R}\dif t\right)\\
&\leq\mE\left(\int_0^{T}|f(t,X^R_t(x))|^\theta 1_{|X^R_t(x)|<R}\dif t\right)
\leq C_R\left(\int_0^{T}\!\!\!\int_{B_R}|f(t,y)|^q\dif y\dif t\right)^{\frac{\theta}{q}}.
\end{align*}
Thus, using the same technique as in the proof of Lemma \ref{Le33} and by (\ref{y2}), we get \eqref{k22}.
\end{proof}

The following lemma will be used in the proof of irreducibility.

\bl\label{l3} For given $x_0, y_0\in\mR^d$ and $m\geq 1$, let $Y_t$ solve the following SDE:
\begin{align}
\label{Eq5}
\dif Y_t=-m(Y_t-y_0)\dif t+b(t,Y_t)\dif t+\sigma(t,Y_t)\dif W_t,~~Y_0=x_0.
\end{align}
Under {\bf(H1)} and {\bf(H2)}, for any $T>0$, there exist constants $C_0, C_1>0$ such that for all $t\in[0,T]$ and $m\geq 1$,
\begin{align}
\mE|Y_t-y_0|^2\leq C_0\e^{-mt}|x_0-y_0|^2+\frac{C_1}{\sqrt{m}}\label{EB4}
\end{align}
and
\begin{align}
\mE\left(\sup_{t\in[0,T]}|Y_t|^2\right)<\infty.\label{EB5}
\end{align}
\el
\begin{proof}
Let $\tilde b(t,x):=-m(x-y_0)+b(t,x)$ and $f$ be as in (\ref{EB2}). As in the calculations of (\ref{EB1}), we have
\begin{align*}
(\p_s f+\sL^{\sigma,\tilde b}_s f)(s,x)
&\leq\alpha(1+|x|^2)^\alpha\e^{-\lambda s}f(s,x)\left(-\tfrac{\lambda}{\alpha}+\tfrac{2\<\tilde b,x\>+2(\alpha+1)(1+|x|^2)^\alpha\|\sigma\|^2}{1+|x|^2}\right)\\
&\!\!\stackrel{(\ref{c1})}{\leq}\alpha(1+|x|^2)^\alpha\e^{-\lambda s}f(s,x)\left(-\tfrac{\lambda}{\alpha}+2C_{\alpha+1}+\tfrac{2m(|x|\cdot|y_0|-|x|^2)}{1+|x|^2}\right).
\end{align*}
If $|x|\leq |y_0|$, then
$$
(\p_s f+\sL^{\sigma,\tilde b}_s f)(s,x)\leq\alpha(1+|y_0|^2)^\alpha\e^{-\lambda s}\exp\Big(\e^{-\lambda s}(1+|y_0|^2)^{\alpha}\Big)
\left\{2C_{\alpha+1}+2m|y_0|^2\right\}.
$$
If $|x|>|y_0|$ and choose $\lambda>2\alpha C_{\alpha+1}$, then
$$
(\p_s f+\sL^{\sigma,\tilde b}_s f)(s,x)\leq 0.
$$
Hence,
\begin{align}
\mE \exp\Big\{\e^{-\lambda t}(1+|Y_t|^2)^{\alpha}\Big\}=\mE f(t,Y_t)\leq f(0,x_0)+C(y_0)(1+m)t.\label{EB3}
\end{align}
On the other hand,
by It\^o's formula, we have for all $t\in[0,T]$,
\begin{align*}
\frac{\dif\mE|Y_t-y_0|^2}{\dif t}&=-2m\mE|Y_t-y_0|^2+2\mE\<Y_t-y_0, b(t,Y_t)\>+\mE\|\sigma(t,Y_t)\|^2\\
&\stackrel{(\ref{c1})}{\leq} -2m\mE|Y_t-y_0|^2+2C_{1/2}(1+\mE|Y_t|^2)+2|y_0|\mE|b(t,Y_t)|\\
&\stackrel{(\ref{c3})}{\leq} 2(C_{1/2}-m)\mE|Y_t-y_0|^2+C+\left(\mE\exp\{C(1+|Y_t|^2)^{\alpha'}\}\right)^{\frac{1}{2}}\\
&\stackrel{(\ref{EB3})}{\leq} 2(C_{1/2}-m)\mE|Y_t-y_0|^2+C\sqrt{m},
\end{align*}
where $C=C(T, x_0,y_0)$ is independent of $m$.
By Gronwall's inequality, we have
$$
\mE|Y_t-y_0|^2\leq \e^{2(C_2-m)t}|x_0-y_0|^2+C\sqrt{m}\e^{2(C_2-m)t}\int^t_0 \e^{2(m-C_2)s}\dif s,
$$
which then gives (\ref{EB4}).
As for (\ref{EB5}), it follows by (\ref{Eq5}), (\ref{c3}) and (\ref{EB3}).
\end{proof}

We are now in a position to give

\begin{proof}[Proof of Theorem \ref{main1}]
For any $p\geq 2$ and $T>0$, by (\ref{c3}) and Burkholder's inequality, we have for all $0\leq s\leq t\leq T$,
\begin{align}\label{RG4}
\begin{split}
&\mE|X_t(x)-X_s(x)|^p\\
&\quad\leq C\mE\left(\int^t_s |b(r,X_r(x))|\dif r\right)^p+C\mE\left|\int^t_s \sigma(r,X_r(x))\dif W_r\right|^p\\
&\quad\leq C\mE\left(\int^t_s |b(r,X_r(x))|^2\dif r\right)^{\frac{p}{2}}+C\mE\left(\int^t_s\|\sigma(r,X_r(x))\|^2\dif r\right)^{\frac{p}{2}}\\
&\quad\leq C(t-s)^{\frac{p}{2}-1}\mE\int^t_s(|b(r,X_r(x))|^p+\|\sigma(r,X_r(x))\|^p)\dif r\\
&\quad\leq C(t-s)^{\frac{p}{2}-1}\mE\int^t_s\Big(1_{\alpha>0}\exp\Big\{C(1+|X_r(x)|^2)^{\alpha'}\Big\}\\
&\qquad\qquad\qquad\qquad\qquad+1_{\alpha=0}(1+|X_r(x)|^2)^{p\gamma_2}\Big)\dif r\\
&\quad\leq C(t-s)^{\frac{p}{2}}\left(1_{\alpha>0}\exp\Big\{(1+|x|^2)^{\alpha}\Big\}+1_{\alpha=0}(1+|x|^2)^{p\gamma_2}\right),
\end{split}
\end{align}
where the last step is due to $\alpha'\in[0,\alpha)$, Young's inequality and Lemma \ref{Le31}.

Next, set $Z_t:=X_t(x)-X_t(y)$. For any $p\geq 1$,  by It\^o's formula we have
\begin{align}
\begin{split}
|Z_t|^{2p}&=|x-y|^{2p}+2p\int_0^t|Z_s|^{2(p-1)}\Big\langle Z_s,\big[\sigma(s,X_s(x))-\sigma(s,X_s(y))\big]\dif W_s\Big\rangle\\
&\quad+2p\int_0^t|Z_s|^{2(p-1)}\Big\langle Z_s,\big[b(s,X_s(x))-b(s,X_s(y))\big]\Big\rangle\dif s\\
&\quad+2p\int_0^t|Z_s|^{2(p-1)}\|\sigma(s,X_s(x))-\sigma(s,X_s(y))\|^2\dif s\\
&\quad+2p(p-1)\int_0^t|Z_s|^{2(p-2)}\big|\big[\sigma(s,X_s(x))-\sigma(s,X_s(y))\big]^*Z_s\big|^2\dif s\\
&=:|x-y|^{2p}+\int_0^t|Z_s|^{2p}\Big(\xi(s)\dif W_s+\eta(s)\dif s\Big),
\end{split}\label{eq1}
\end{align}
where
$$
\xi(s):=\frac{2p\big[\sigma(s,X_s(x))-\sigma(s,X_s(y))\big]^*Z_s}{|Z_s|^2}
$$
and
\begin{align*}
\eta(s):&=\frac{2p\big\langle Z_s,[b(s,X_s(x))-b(s,X_s(y))]\big\rangle}{|Z_s|^2}\\
&+\frac{2p\|\sigma(s,X_s(x))-\sigma(s,X_s(y))\|^2}{|Z_s|^2}\\
&+\frac{2p(p-1)|[\sigma(s,X_s(x))-\sigma(s,X_s(y))]^*Z_s|^2}{|Z_s|^4}.
\end{align*}
Here we use the convention $\frac{0}{0}=0$.
By Dol\'eans-Dade's exponential formula, we have
$$
|Z_t|^{2p}=|x-y|^{2p}\exp\Bigg\{\int_0^t\xi(s)\dif W_s-\frac{1}{2}\int_0^t|\xi(s)|^2\dif s+\int_0^t\eta(s)\dif s\Bigg\}.
$$
For $R>|x|\vee|y|$, define a stopping time
$$
\tau_R:=\inf\Big\{t\geq 0: |X_t(x)|\vee|X_t(y)|\geq R\Big\}.
$$
By \eqref{UY5}, we have for $s<\tau_R$,
$$
|\xi(s)|\leq 2^{d+1}p\Big(\cM_{2R}|\nabla\sigma(t,\cdot)|(X_s(x))+\cM_{2R}|\nabla\sigma(t,\cdot)|(X_s(y))\Big).
$$
Since $\cM_{2R}|\nabla\sigma(t,\cdot)|(x)\in L^q_{loc}(\mR_+\times\mR^d)$ with $q>d+2$, by Lemma \ref{Le34}, we have for any $\kappa>0$,
$$
\mE\exp\Bigg\{\kappa\int_0^{T\wedge\tau_R}|\xi(s)|^2\dif s\Bigg\}<\infty.
$$
Hence, for any $\kappa>0$, by Novikov's criterion,
$$
t\mapsto\exp\Bigg\{\kappa\int_0^{t\wedge\tau_R}\xi(s)\dif W_s-\frac{\kappa^2}{2}\int_0^{t\wedge\tau_R}|\xi(s)|^2\dif s\Bigg\}=:\sE_\kappa(t)
$$
is a martingale. Thus, by H\"older's inequality we have
\begin{align}
\begin{split}
\mE|Z_{t\wedge\tau_R}|^{2p}&\leq|x-y|^{2p}\Big(\mE\sE_2(t)\Big)^{\frac{1}{2}}\Bigg(\mE\exp\left\{
\int_0^{t\wedge\tau_R}\Big(|\xi(s)|^2+2\eta(s)\Big)\dif s\right\}\Bigg)^{\frac{1}{2}}\\
&=|x-y|^{2p}\Bigg(\mE\exp\left\{\int_0^{t\wedge\tau_R}\Big(|\xi(s)|^2+2\eta(s)\Big)\dif s\right\}\Bigg)^{\frac{1}{2}}.
\end{split}\label{UY9}
\end{align}
On the other hand, in view of
\begin{align*}
|Z_s|^2(|\xi(s)|^2+2\eta(s))&\leq 8p^2\|\sigma(s,X_s(x))-\sigma(s,X_s(y))\|^2\\
&\quad+4p\big\langle Z_s,[b(s,X_s(x))-b(s,X_s(y))]\big\rangle\\
&\!\!\stackrel{(\ref{c4})}{\leq} 4p|Z_s|^2(F_{2p}(s,X_s(x))+F_{2p}(s,X_s(y))),
\end{align*}
by \eqref{UY9}, (\ref{cf}) and Lemma \ref{Le33}, as well as Lemma \ref{Le31} and Fatou's lemma, we further have
\begin{align}
\begin{split}
\mE|Z_t|^{2p}
&\leq|x-y|^{2p}\Bigg(\mE\exp\left\{4p\int_0^{t}\Big(F_{2p}\big(s,X_s(x)\big)+F_{2p}\big(s,X_s(y)\big)\Big)\dif s\right\}\Bigg)^{\frac{1}{2}}\\
&\leq C|x-y|^{2p}\Big\{g(x)g(y)\Big\}^{\frac{1}{4}}\leq C|x-y|^{2p}\Big(g(x)^{\frac{1}{2}}+g(y)^{\frac{1}{2}}\Big),
\end{split}\label{RG3}
\end{align}
where $g(x):=1_{\alpha>0}\e^{(1+|x|^2)^\alpha}+1_{\alpha=0}(1+|x|^2)^\gamma$,
which, together with  (\ref{RG4}) and Kolmogorov's continuity criterion, yields that $X_t(x)$ admits a bicontinuous version, and for any $T,R>0$ and $p\geq 1$,
\begin{align}
\mE\left(\sup_{t\in[0,T],|x|\leq R}|X_t(x)|^p\right)<+\infty.\label{EQ09}
\end{align}
\\
{\bf (A)} It follows by (\ref{RG3}) and Lemma \ref{Le21} with $q=\infty$,  $U$ being any ball.
\\
\\
{\bf (B)} Following the above proof, for any $T>0$, by H\"older's inequality and Doob's maximal inequality, we have
\begin{align*}
\mE\left(\sup_{t\in[0,T\wedge\tau_R]}|Z_t|^{2p}\right)&\leq|x-y|^{2p}\left(\mE\sup_{t\in[0,T]}
\sE^2_1(t)\right)^{\frac{1}{2}}\Bigg(\mE\exp\left\{\sup_{t\in[0,T\wedge\tau_R]}2\int_0^t\eta(s)\dif s\right\}\Bigg)^{\frac{1}{2}}\\
&\leq 2|x-y|^{2p}\left(\mE\sE^2_1(T)\right)^{\frac{1}{2}}\Bigg(\mE\exp\left\{\sup_{t\in[0,T\wedge\tau_R]}2\int_0^t\eta(s)\dif s\right\}\Bigg)^{\frac{1}{2}}\\
&\leq 2|x-y|^{2p}\Big(\mE\sE_{4}(T)\Big)^{1/4}\Bigg(\mE\exp\left\{6\int_0^{T}|\xi(s)|^2\dif s\right\}\Bigg)^{1/4}\\
&\quad\times\Bigg(\mE\exp\left\{\sup_{t\in[0,T\wedge\tau_R]}2\int_0^t\eta(s)\dif s\right\}\Bigg)^{1/2}.
\end{align*}
By the additional assumption \eqref{Sigma}, as in the above proof, we get
$$
\mE\left(\sup_{t\in[0,T]}|Z_t|^{2p}\right)
\leq C|x-y|^{2p}\Big\{g(x)g(y)\Big\}^{\frac{1}{4}}\leq C|x-y|^{2p}\Big(g(x)^{\frac{1}{2}}+g(y)^{\frac{1}{2}}\Big).
$$
which, together with Lemma \ref{Le21} with $q,r=\infty$ and  $U$ being any ball, implies \eqref{so1}.
\\
\\
{\bf (C)} For each $n\in\mN$, let $\chi_n(x)$ be a nonnegative smooth function in $\mR^d$ with $\chi_n(x)=1$ for all $x\in B_n$
and $\chi_n(x)=0$ for all $x\notin B_{n+1}$. Let
$$
b_n(t,x):=\chi_n(x)b(t,x),\ \ \sigma_n(t,x):=\sigma(t,\chi_n(x)x).
$$
Clearly, for any $T>0$,
$$
b_n\in L^{q}([0,T]\times\mR^d), \ \ \nabla\sigma_n\in L^{q}([0,T]\times\mR^d),
$$
and for some $K_n>0$,
$$
K^{-1}_n|\xi|\leq|\sigma_n(t,x)\xi|\leq K_n|\xi|,\ \ (t,x)\in [0,T]\times\mR^d,\ \ \xi\in\mR^m.
$$
Let $X^n_t(x)$ be the solution of SDE (\ref{1}) corresponding to $b_n$ and $\sigma_n$.
By \cite[Theorem 1.1]{Zh1} or \cite{Zh4}, for any bounded measurable function $f$ and $t>0$,
\begin{align}
x\mapsto \mE f(X^n_t(x))\mbox{ is continuous.}\label{EB6}
\end{align}
Fix $R>0$. For $n>R$, define a stopping time
$$
\tau_{n,R}:=\left\{t\geq 0: \sup_{|x|\leq R}|X_t(x)|\geq n\right\}.
$$
By Chebyshev's inequality and (\ref{EQ09}), we have
\begin{align}
\lim_{n\to\infty}\mP(t>\tau_{n,R})\leq\lim_{n\to\infty}\mE\left(\sup_{s\in[0,t],|x|\leq R}|X_s(x)|^p\right)/n=0.\label{EB7}
\end{align}
Moreover, by the local uniqueness of solutions to SDE (\ref{1}) (see \cite{Zh1}), we have
$$
X_t(x)=X^n_t(x),\ \ |x|\leq R,\ \ t\in[0,\tau_{n,R}].
$$
Let $f$ be a bounded measurable function. For any $x,y\in B_R$, we have
\begin{align*}
|\mE(f(X_t(x))-f(X_t(y)))|&\leq\big|\mE\big(f(X_t(x))-f(X_t(y))1_{t\leq\tau_{n,R}}\big)\big|+2\|f\|_\infty \mP(t>\tau_{n,R})\\
&=\big|\mE\big(f(X^n_t(x))-f(X^n_t(y))1_{t\leq\tau_{n,R}}\big)\big|+2\|f\|_\infty \mP(t>\tau_{n,R})\\
&\leq\big|\mE\big(f(X^n_t(x))-f(X^n_t(y))\big)\big|+4\|f\|_\infty \mP(t>\tau_{n,R}),
\end{align*}
which together with (\ref{EB6}) and (\ref{EB7}) yields the continuity of $x\mapsto \mE(f(X_t(x)))$.
\\
\\
{\bf (D)} Our proof is adapted from \cite{Re-Wu-Zh}. It suffices to prove that for any $T, a>0$ and $x_0,y_0\in\mR^d$,
$$
\mP(|X_T(x_0)-y_0|\leq a)>0.
$$
In what follows, we shall fix $T, a>0$ and $x_0,y_0\in\mR^d$. Let $Y_t(x_0)$ solve SDE (\ref{Eq5}) and for $N>0$, set
$$
\tau_N:=\inf\{t: |Y_t(x_0)|\geq N\}.
$$
By (\ref{EB4}) and (\ref{EB5}), we may choose $N$ and $m$ large enough so that
\begin{align}
\label{qandp}
\mP(\tau_N\leq T)+\mP(|Y_T(x_0)-y_0|>a)<1.
\end{align}
Define
$$
U_t:=-m~\sigma(t,Y_t)^{*}[\sigma(t,Y_t)\sigma(t,Y_t)^{*}]^{-1}(Y_t-y_0)
$$
and
$$Z_T:=\exp\left(\int_0^{T\wedge\tau_N}U_s\dif W_s-\frac{1}{2}\int_0^{T\wedge\tau_N}|U_s|^2\dif s\right).
$$
Since $|U_{t\wedge \tau_N}|^2$ is bounded, $\mE[Z_T]=1$ by Novikov's criteria.

By Girsanov's theorem, $\tilde W_t:=W_t+V_t$ is a $\mQ$-Brownian motion, where
$$
V_t:=\int_0^{t\wedge\tau_N}U_s\dif s,\quad \mQ:=Z_T\mP.
$$
By (\ref{qandp}) we have
\begin{align}
\label{qless1} \mQ(\{\tau_N\leq T\}\cup\{|Y_T(x_0)-y_0|>a\})<1.
\end{align}
Notice that the solution $Y_t$ of (\ref{Eq5}) also solves the following SDE:
\begin{align*}
Y_{t\wedge\tau_N}=x_0+\int_0^{t\wedge \tau_N}b(s,Y_s)\dif s+\int_0^{t\wedge\tau_N}\sigma(s,Y_s)\dif \tilde W_s.
\end{align*}
Set
$$
\theta_N:=\inf\{t: |X_t|\geq N\}.
$$
Then the uniqueness in distribution for (\ref{1}) yields that the law of $\{(X_t1_{\{\theta_N\geq t\}})_{t\in [0,T]},\theta_N\}$ under $\mP$ is the same as that of
$\{(Y_t1_{\{\tau_N\geq t\}})_{t\in [0,T]},\tau_N\}$ under $\mQ$. Hence
\begin{align*}
\mP(|X_T(x_0)-y_0|>a)&\leq\mP(\{\theta_N\leq T\}\cup\{\theta_N\geq T, |X_T(x_0)-y_0|>a\})\\
&=\mQ(\{\tau_N\leq T\}\cup\{\tau_N\geq T, |Y_T(x_0)-y_0|>a\})\\
&\leq \mQ(\{\tau_N\leq T\}\cup\{|Y_T(x_0)-y_0|>a\})<1.
\end{align*}
The proof is complete.
\end{proof}

\section{Proof of Theorem \ref{main3}}


We first prepare the following easy lemma.
\bl\label{Le41}
Let $f:\mR^d\to\mR$ be a measurable function. Assume that for some $g_1\in L^{p_1}_{loc}(\mR^d)$,
 $g_2\in L^{p_2}_{loc}(\mR^d)$ and some $R>0$,
\begin{align}
|f(x)-f(y)|&\leq |x-y|(g_1(x)+g_1(y)),\ \ \ \forall x,y\in B_{3R},\\
|f(x)-f(y)|&\leq |x-y|(g_2(x)+g_2(y)),\ \ \ \forall x,y\notin B_R.
\end{align}
Then we have for all $x,y\in\mR^d$ with $|x-y|\leq R$,
$$
|f(x)-f(y)|\leq 2^{d+1}|x-y|(g(x)+g(y)),
$$
where
$$
g(x)=\cM_R g_1(x) 1_{|x|\leq 2R}+\cM_R g_2(x) 1_{|x|>2R},
$$
and $\cM_R g_i(x):=\sup_{s\in(0,R)}\fint_{B_s}|g_i(x+z)|\dif z$, $i=1,2$.
\el
\begin{proof}
First of all, by the assumptions and Lemma \ref{Le21}, we have
$$
|\nabla f(x)|\leq 2g_1(x),\ |x|<3R, \ |\nabla f(x)|\leq 2g_2(x),\ |x|>R.
$$
By \eqref{UY5}, we have for Lebesgue-almost all $x,y\in\mR^d$ with $|x-y|<R$,
$$
|f(x)-f(y)|\leq 2^d|x-y|\left(\cM_R|\nabla f|(x)+\cM_R|\nabla f|(y)\right),
$$
which in turn implies the desired estimate by the definition of $\cM_R$ and redefinition of $g(x)$ on a Lebesgue zero set.
\end{proof}

Below, we fix $T>0$ and write for $p\in[1,\infty]$,
$$
\mL^p(T):=L^p([0,T]\times\mR^d).
$$
Let $\chi\in C^\infty(\mR^d;[0,1])$ be a cutoff function with
$$
\chi(x)=1,\forall |x|\leq 1,\ \ \chi(x)=0, \forall |x|>2,\ \ \|\nabla\chi\|_\infty\leq 2,
$$
and for $R>0$, we set
$$
\chi_R(x):=\chi(x/R),\ \ \bar\chi_R(x)=1-\chi_R(x).
$$
Let $R_0$ be as in ${\bf {(H2')}}$. Without loss of generality, we may assume $R_0\geq 4$ so that
\begin{align}
\|\nabla\chi_{R_0}\|_\infty\leq\|\nabla\chi\|_\infty/R_0\leq 1/2.\label{EQ00}
\end{align}
We make the following decomposition for $b$:
$$
b=b_1+b_2,\ \ b_1:=b\chi_{R_0},\ \ b_2:=b\bar\chi_{R_0}.
$$
In view of ${\bf{(H1')}}$, the function $b_1$ is global $L^q$-integrable; while ${\bf{(H2')}}$ implies that $b_2$ satisfies ${\bf{(H2)}}$.
On the other hand, by Sobolev's embedding theorem, ${\bf{(H1')}}$ and ${\bf{(H2')}}$ also imply that for some $\alpha\in(0,1)$ and $C>0$,
$$
\|\sigma(t,x)-\sigma(t,y)\|\leq C|x-y|^\alpha.
$$

The following result is an easy combination of \cite[p.120, Theorem 1]{Kr0} and \cite[Theorem 10.3 and Lemma 10.2]{Kr-Ro}
(see \cite[Theorem 3.5]{Zh4} for a detailed proof).
\bl
Let $q>d+2$. Under ${\bf{(H1')}}$ and ${\bf{(H2')}}$, for any $\lambda>0$, there exists a unique solution $u\in\mL^q(T)$
with $\nabla^2 u\in\mL^q(T)$ to the following backward PDE:
\begin{align}
\p_t u+ \sL^{\sigma,b_1}_t u+b_1=\lambda u,\quad u(T)=0.\label{eq}
\end{align}
Moreover, there exist a $\lambda>0$ and a positive constant $C=C(K,d,q,T,\lambda,\|b_1\|_{\mL^q(T)})$ such that
\begin{align}
\|\p_tu\|_{\mL^q(T)}+\|\nabla^2 u\|_{\mL^q(T)}\leq C<\infty\,\,  \mbox{ and }\,\,\|u\|_{\mL^\infty(T)}+\|\nabla u\|_{\mL^\infty(T)}\leq \tfrac{1}{2}.\label{es1}
\end{align}
\el

Let $u(t,x)$ be as in the above lemma. Now, we want to follow the same idea as in \cite{Zh1} to perform
Zvonkin's transformation and transform SDE (\ref{1}) into a new one with coefficients
satisfying ${\bf{(H1)}}$ and ${\bf{(H2)}}$. However, if we argue entirely the same as usual and consider the transform
$$
(t,x)\mapsto \Psi(t,x):=x+u(t,x),
$$
then one finds that condition (\ref{cf}) may not be satisfied for the new coefficients (see Lemma \ref{Le43} below).
For this reason, we define
$$
u_{R_0}(t,x):=u(t,x)\chi_{2R_0}(x),\quad \Phi_t(x):=x+u_{R_0}(t,x),
$$
where $R_0$ is the same as in ${\bf {(H2')}}$.

We have
\bl\label{Le42}
The following statements hold:
\begin{enumerate}[(i)]
\item For each $t\in[0,T]$, the map $x\rightarrow\Phi_t(x)$ is a $C^1$-diffeomorphism and
$$
\|\nabla\Phi_t\|_\infty,\|\nabla\Phi^{-1}_t\|_\infty\leq 2.
$$
Moreover, $\nabla\Phi_t(x)$ and $\nabla\Phi^{-1}_t(x)$ are H\"older continuous in $x$ uniformly in $t\in[0,T]$.
\item Let $q>d+2$. We have $\p_t\Phi_t,\nabla^2\Phi_t, \p_t\Phi^{-1}_t,\nabla^2\Phi^{-1}_t\in\mL^q(T)$ and
\begin{align}
\p_t\Phi_t+\sL^{\sigma,b_1}_t\Phi_t=\sigma_{ik}\sigma_{jk}\p_i u\p_j\chi_{2R_0}+\tfrac{1}{2}u \sigma_{ik}\sigma_{jk}\p_i\p_j\chi_{2R_0}
+\lambda u_{R_0}.\label{RT1}
\end{align}
Here and below we use Einstein's convention for summation.
\end{enumerate}
\el
\begin{proof}
By (\ref{es1}) and (\ref{EQ00}), we have
$$
\tfrac{1}{2}|x-y|\leq|\Phi_t(x)-\Phi_t(y)|\leq \tfrac{3}{2}|x-y|.
$$
Thus, (i) follows by (\ref{es1}) and Sobolev's embedding result (see \cite[Lemma 10.2]{Kr-Ro}).
\\
\\
(ii) $\p_t\Phi_t,\nabla^2\Phi_t, \p_t\Phi^{-1}_t,\nabla^2\Phi^{-1}_t\in\mL^q(T)$ follows by (\ref{es1}) and (i). Moreover, by elementary calculations, we have
\begin{align}
\begin{split}
\p_t\Phi_t+\sL^{\sigma,b_1}_t\Phi_t&=\sigma_{ik}\sigma_{jk}\p_i u\p_j\chi_{2R_0}+\tfrac{1}{2}u \sigma_{ik}\sigma_{jk}\p_i\p_j\chi_{2R_0}\\
&+ub^i_1\p_i\chi_{2R_0}+\lambda u_{R_0}+b_1(1-\chi_{2R_0}).
\end{split}\label{RT3}
\end{align}
Notice that
$$
b^i_1\p_i\chi_{2R_0}=\chi_{R_0}b^i\p_i\chi_{2R_0}=0,\ \ b_1(1-\chi_{2R_0})=b\chi_{R_0}(1-\chi_{2R_0})=0.
$$
Equality (\ref{RT1}) follows by (\ref{RT3}).
\end{proof}

Basing on the above lemma, we may prove the following Zvonkin transformation (see \cite{F-G-P,Zh1} for more details).
\bl\label{Le43}
Let $h$ be defined by the right hand side of (\ref{RT1}).
Then $X_t$ solves SDE (\ref{1}) if and only if $Y_t:=\Phi_t(X_t)$ solves the following SDE:
\begin{align}
\dif Y_t=\tilde{b}(t,Y_t)\dif t+\tilde{\sigma}(t,Y_t)\dif W_t,    \label{03}
\end{align}
where
\begin{align}
\tilde\sigma:=(\nabla\Phi\cdot\sigma)\circ\Phi^{-1},\ \ \tilde{b}:=(h+b_2\cdot\nabla\Phi)\circ\Phi^{-1}.\label{RT4}
\end{align}
\el
\begin{proof}
($\Rightarrow$) By (\ref{RT1}) and generalized It\^o's formula (see \cite[p.122, Theorem 1]{Kr}), we have (\ref{03}).
\\
\\
($\Leftarrow$)
By elementary calculations, it is easy to check that
$$
\p_t\Phi^{-1}_t+\sL^{\tilde\sigma,\tilde b}_t\Phi^{-1}_t=b\circ\Phi_t^{-1}.
$$
As above, using generalized It\^o's formula again, we obtain that $\Phi^{-1}_t(Y_t)$ solves SDE (\ref{1}).
\end{proof}

Now we give

\begin{proof}[Proof of Theorem \ref{main3}]
Since {\bf (A)}-{\bf (D)} are invariant under diffeomorphism transformation $x\mapsto\Phi_t(x)$,
by Lemma \ref{Le43}, it suffices to check that $\tilde\sigma$ and $\tilde b$ defined by (\ref{RT4})
satisfy {\bf (H1)}-{\bf (H2)} so that we can use Theorem \ref{main1} to conclude the proof.

First of all, {\bf (H1)} is obvious by Lemma \ref{Le42} and ${\bf (H1')}$.
For {\bf (H2)}, by definitions \eqref{RT4} and Lemma \ref{Le42}, it is easy to see that
\begin{align*}
|\tilde b(t,x)|\leq \|h\|_\infty+2|b_2(t,\Phi^{-1}_t(x))|\leq C(1+|\Phi^{-1}_t(x)|)\leq C(1+|x|),
\end{align*}
and by (\ref{Mo2}), for any $R>0$, there are functions $g_R, \hat g_R\in L^q_{loc}(\mR_+\times\mR^d)$ such that for all $x,y\in \mR^d$ with $|x-y|\leq R$
\begin{align*}
\|\tilde\sigma(t,x)-\tilde\sigma(t,y)\|&\leq |x-y|(g_R(t,x)+g_R(t,y)),\\
|\tilde b(t,x)-\tilde b(t,y)|&\leq |x-y|(\hat g_R(t,x)+\hat g_R(t,y)).
\end{align*}
On the other hand, by the definition of $\Phi$, there exists a $R_1\geq 2R_0$ large enough such that
$$
\Phi_t(x)=\Phi^{-1}_t(x)=x,\ \ |x|\geq R_1.
$$
Hence, for $|x|\geq R_1$, we have
$$
\tilde b(t,x)=b(t,x),\ \ \tilde\sigma(t,x)=\sigma(t,x).
$$
Thus, by  {\bf (H2$'$)} and Lemma \ref{Le41}, one sees that {\bf (H2)} and \eqref{Sigma} hold for $\tilde b$ and $\tilde\sigma$.
The proof is complete.
\end{proof}

{\bf Acknowledgements.} The authors would like to express their deep thanks to the referee for carefully reading the manuscript and useful suggestions. This work is supported by NNSFs of China (Nos. 11271294, 11325105) and the Fundamental Research Funds
for the Central Universities (No. 2014201020208).

\end{document}